\documentclass{article}

\usepackage{amsthm,amsmath,amssymb}

\newtheorem{theorem}{Theorem}[section]

\newtheorem{proposition}[theorem]{Proposition}

\newtheorem{remark}{Remark} [section]
\newtheorem{corrollary}[theorem]{Corollary}

\DeclareMathOperator{\Hess}{Hess}

\DeclareMathOperator{\Var}{Var}

\DeclareMathOperator{\id}{Id}
\DeclareMathOperator{\Tr}{Tr}

\newcommand{\R}{\mathbb{R}}

\usepackage{color}

\begin{document}

\title{Stein kernels and moment maps}

\author{Max Fathi\thanks{CNRS and Institut de Math\'ematiques de Toulouse, Universit\'e de Toulouse, max.fathi@math.univ-toulouse.fr}}

\maketitle

\begin{abstract}
We describe a construction of Stein kernels using moment maps, which are solutions to a variant of the Monge-Amp\`ere equation. As a consequence, we show how regularity bounds in certain weighted Sobolev spaces on these maps control the rate of convergence in the classical central limit theorem, and derive new rates in Kantorovitch-Wasserstein distance in the log-concave situation, with explicit polynomial dependence on the dimension. 
\end{abstract}

\section{Introduction}
Stein's method is a set of techniques introduced by Stein \cite{Ste72, Ste86} to estimate distances between probability measures. We refer to the survey \cite{Ros11} for an overview. We shall be interested in one particular way of implementing Stein's method in the Gaussian setting, based on the notion of Stein kernels. Let $\mu$ be a probability measure on $\R^d$. A matrix-valued function $\tau_{\mu} : \R^d \longrightarrow \mathcal{M}_d(\R)$ is said to be a Stein kernel for $\mu$ (with respect to the standard Gaussian measure $\gamma$ on $\R^d$) if for any smooth test function $f$ taking values in $\mathbb{R}^d$,  we have
\begin{equation} \label{eq_stein}
\int{x \cdot f d\nu} = \int{\langle \tau_{\mu}, \nabla f \rangle_{\mathrm{HS}} d\nu}. 
\end{equation}
For  applications, it is generally enough to consider the restricted class of test functions $f$ satisfying $\int (|f|^2 + \|\nabla f\|_{\mathrm{HS}}^2 )d\mu <\infty$, in which case both integrals in \eqref{eq_stein} are well-defined as soon as  $\tau_{\mu} \in L^2(\mu)$, provided $\mu$ has finite second moments.  

The  motivation behind the definition is that, since the Gaussian measure is the only probability distribution satisfying the integration by parts formula
\begin{equation} \label{ibp_gauss}
\int{x \cdot f d\gamma} = \int{\operatorname{div} (f)d\gamma},
\end{equation} 
a Stein kernel $\tau_{\mu}$ coincides with the identity matrix, denoted by $\id$, if and only if the measure $\mu$ is  equal to $\gamma$. Hence,  the Stein kernel can be used to control how far $\mu$ is from being a standard Gaussian measure in terms of how much it violates the integration by parts formula \eqref{ibp_gauss}. It appears implicitly in many works on Stein's method, and has recently been the topic of more direct investigations \cite{AMV10, NPR10, NPS14a, LNP15, LVY, CFP17}. The one-dimensional case, where Stein kernels can be explicitly constructed from the density, has been extensively studied \cite{LRS17}. It has applications to central limit theorems \cite{NP12}, concentration inequalities \cite{NV09, LNP15, Sau18} and random matrix theory \cite{Cha09}. 

A related quantity is the Stein discrepancy
$$S(\mu)^2 := \underset{\tau}{\inf} \int{|\tau - \id|^2d\mu}$$
where the infimum is taken over all possible Stein kernels for $\mu$, since they may not be unique. This quantity has two main interesting properties: it controls the $L^2$ Kantorovitch-Wasserstein distance to the Gaussian \cite{LNP15}, and is monotone along the central limit theorem \cite{CFP17}. 

The aim of this work is to describe how we can construct Stein kernels using a correspondence between centered measures and convex functions, known as the moment measure problem, or moment map problem, which we shall describe in Section \ref{sect_moment_maps}. The main motivation was to give a construction of Stein kernels using optimal transport maps, of which these moment maps can be viewed as a variant. The Stein kernels we shall build have several nice properties that do not seem to be necessarily satisfied by previous constructions. Most notably they shall always takes values that are symmetric, nonnegative matrices. As an application, we shall derive in Section \ref{sect_clt} new bounds on the rate of convergence in the multi-dimensional central limit theorem when the random variables are log-concave, with explicit dependence on the dimension. Their main interest is that the dependence on the dimension will improve on a more general result of Bonis \cite{Bon16} in the particular case of log-concave measures. In particular, we shall derive the sharp dependence on the dimension in the uniformly log-concave setting. In section 4, we shall discuss a multi-dimensional generalization of a result of Saumard \cite{Sau18} on weighted Poincar\'e inequalities involving Stein kernels, and in Section 5 we shall briefly point out a construction of Stein kernels with respect to non-Gaussian reference measures. 

\section{Stein kernels and moment maps}

\subsection{Moment maps} \label{sect_moment_maps}

In \cite{CK14} (revisited in \cite{San16}, and following earlier works \cite{WZ04, Don08, BB13, Leg16}), the following theorem was established: 

\begin{theorem}[Cordero-Erausquin and Klartag 2015] \label{thm_moment_measures}
Let $\mu$ be a centered measure, with finite first moment and that is not supported on a hyperplane. Then there exists a convex function $\varphi$ such that $\mu$ is the pushforward of the probability measure with density $e^{-\varphi}$ by the map $\nabla \varphi$. This function $\varphi$ is called the moment map of $\mu$.
\end{theorem}

This result can be seen as a variant of the optimal transport problem, where instead of specifying two measures, we fix a target measure, and look for both an original measure and a transport map while imposing the constraint that the map should be the gradient of the potential of the measure. Indeed, here $\nabla \varphi$ is also the Brenier map from optimal transport theory \cite{Vil03} sending $e^{-\varphi}$ onto $\mu$. 

The convex function given by this theorem may well not be smooth, most notably when $\mu$ is a combination of Dirac masses. For example, if $\mu$ is the uniform measure on $\{-1, +1\}$, viewed as a subset of $\R$, the convex function is $\varphi(x) = |x|$ on $\R$, which is not smooth at the origin. This will cause some issues later on. We can however assume it satisfies some weak continuity property on the boundary of its support (the notion of essential continuity, which is described in \cite{CK14}). A smooth version of this theorem, under extra assumptions, was previously obtained by Berman and Berndtson \cite{BB13}, with earlier results due to Wang and Zhu \cite{WZ04} and Donaldson \cite{Don08}: 

\begin{theorem}[Berman and Berndtson 2013] \label{thm_moment_measures_smooth}
Assume that $\mu$ is supported on a compact, open convex set, and that it has a smooth density $\rho$ on its support. Assume moreover that $C \geq \rho \geq C^{-1}$ on the whole support, for some positive constant $C$. Then the convex function $\varphi$ of Theorem \ref{thm_moment_measures} is smooth and supported on the whole space $\R^d$. 
\end{theorem}

In this result (which is based on Caffarelli's regularity theory for Monge-Amp\`ere PDEs), the convexity of the support plays an essential role to guarantee smoothness of the map. 

We can reformulate those statements as pertaining to solutions of the PDE
\begin{equation} \label{tkepde}
e^{-\varphi} = \rho(\nabla \varphi)\det( \nabla^2 \varphi ).
\end{equation}

This PDE is a variant of the Monge-Amp\`ere equation, sometimes called the toric K\"ahler-Einstein equation. It has been studied in complex geometry, where it is related to the construction of differential structures with specific properties on toric varieties (i.e. quotients of the complex space $(\mathbb{C}^*)^n$). More recently, it has been studied in \cite{K14, KK17, KK18, KM16}, where it was used to establish functional inequalities for log-concave measures. 

A relevant remark to the connection with Stein's method that we shall describe in the next section is that the standard Gaussian measure is the only fixed point of the map $\mu \rightarrow e^{-\varphi}$, where $\varphi$ is the moment map of $\mu$. So in some sense the moment map already contains some information on how far the measure is from being Gaussian. 

In general, unless the dimension is $1$, solutions to \eqref{tkepde} are not explicit. One particular case where it can be determined is for the uniform measure on the unit cube $[-1, 1]^d$, where the moment map is of the form $\varphi(x) = \underset{i=1}{\stackrel{d}{\sum}} \hspace{1mm} 2\log \cosh(x_i/2) + C$. This can be generalized to uniform measures on centered parallelipipeds by composing this function with the appropriate linear map.

\subsection{The connection with Stein kernels}

For now, assume that $\mu$ has a density with respect to the Lebesgue measure which is strictly positive on its support, and is such that the convex function $\varphi$ given by Theorem \ref{thm_moment_measures} is $C^2$. There exists an optimal transport map sending $\mu$ onto $e^{-\varphi}$, which is necessarily $\nabla \varphi^*$, where $\varphi^*$ is the Legendre transform of $\varphi$. $\varphi^*$ is then also $C^2$: since $\nabla \varphi^*$ is the inverse of $\nabla \varphi$ (this is a property of the Legendre transform) and $\Hess \varphi$ is strictly positive on the whole space, $\nabla \varphi^*$ inherits $C^1$ regularity from $\nabla \varphi$. 

\begin{theorem} \label{main_thm_stein}
If $\mu$ has a density $\rho$ with respect to the Lebesgue measure, and the solution $\varphi$ to the PDE \eqref{tkepde} is $C^2$ and supported on the whole space $\R^d$, then $\Hess  \varphi(\nabla \varphi^*)$ is a Stein kernel for $\mu$. Moreover, the Stein discrepancy satisfies
$$S(\mu)^2 \leq \int{|\Hess \varphi - \id|_{HS}^2e^{-\varphi}dx}.$$

In particular, if $\mu$ is supported on a compact, convex set and has density bounded from above and below by positive constants, this result applies. 
\end{theorem}

The regularity assumptions can be weakened, indeed if $\int{|\Hess \varphi - \id|_{HS}^2e^{-\varphi}dx}$ is finite and $\mu$ has a continuous density, then the result will still hold. For general measures with density and full support, the moment map is only in $W^{2,1}_{loc}$ in the interior of its support \cite{DPF13}, which is not enough to make the proof work. But this is not surprising, since for heavy-tailed random variables the CLT may fail, and this would rule out existence of a Stein kernel belonging to $L^2(\mu)$. For background on regularity theory for Monge-Amp\`ere PDEs, we refer to the lecture notes \cite{Fig17}. 

\begin{remark}
An interesting byproduct of this result is that the Stein kernel obtained in this way takes values that are symmetric and positive matrices. In particular, this explains why the explicit formula for Stein kernels in dimension one defines a nonnegative function. This remark will play an important role later on, notably in Section 4. 
\end{remark}

\begin{remark}
The Stein kernel constructed this way seems to be in general different from the one constructed in \cite{CFP17}. Since when the density is supported on a compact, convex set and has density bounded from above and below by positive constants a Poincar\'e inequality holds, existence of a Stein kernel in that situation was already proven in \cite{CFP17}. It is the particular structure of the kernel we obtain here that makes it interesting, as we will see when obtaining new rates of convergence in the CLT. 
\end{remark}

\begin{proof}
Since $\varphi$ is smooth, we have the Stein equation 
$$\int{\nabla \varphi \cdot f e^{-\varphi}dx} = \int{\Tr(\nabla f)e^{-\varphi}dx}.$$
There is no boundary term remaining when integrating by parts because $\varphi$ grows at least linearly at infinity, since it is convex and $\int{e^{-\varphi}dx} < \infty$ (see for example Lemma 2.1 in \cite{Kla07}). 

Fix $g$ a smooth function, and take $f(x) = g(\nabla \varphi (x))$ in the above equation. We get
$$\int{\nabla \varphi(x) \cdot g(\nabla \varphi(x))e^{-\varphi}dx} = \int{\langle \Hess \varphi, \nabla g(\nabla \varphi) \rangle e^{-\varphi} dx}.$$
Applying the change of variable $y = \nabla \varphi^*(x)$, which sends $\mu$ onto $e^{-\varphi}$, we obtain
$$\int{x \cdot g(x) d\mu} = \int{\langle \Hess  \varphi(\nabla \varphi^*), \nabla g \rangle d\mu}$$
which ensures that $\Hess  \varphi(\nabla \varphi^*) = (\Hess \varphi^*)^{-1}$ is indeed a Stein kernel for $\mu$. 

The bound on the Stein discrepancy is an immediate consequence of the change of variable: since $\Hess  \varphi(\nabla \varphi^*)$ is a Stein kernel, by definition of the Stein discrepancy we have
$$S(\mu)^2 \leq \int{|\Hess  \varphi(\nabla \varphi^*) - \id|_{HS}^2d\mu} = \int{|\Hess \varphi - \id|_{HS}^2e^{-\varphi}dx}.$$
\end{proof}

The well-known Caffarelli contraction theorem \cite{Caf00} states that the Brenier map sending the standard Gaussian map onto a uniformly log-concave measure is Lipschitz. Klartag \cite{K14} proved an analogous estimate for moment maps, which leads to the following bound on Stein kernels in that setting: 

\begin{corrollary} \label{cor_est_stein_unif}
Assume that $\mu$ is uniformly convex, that is it is of the form $e^{-V}dx$ with $\Hess V \geq \epsilon \operatorname{Id}$ for some $\epsilon > 0$. Then there exists a Stein kernel with values that are positive symmetric matrices, and  which is uniformly bounded, that is $||\tau||_{op} \leq \epsilon^{-1}$.
\end{corrollary}

In dimension one, this result was pointed out in \cite{Sau18}. Such pointwise estimates can be used to derive properties of the density and concentration inequalities \cite{NV09} and isoperimetric inequalities \cite{Sau18}.  

\begin{proof}
The Stein kernel described in this statement is the one built in Theorem \ref{main_thm_stein}, all that we need to do is to prove the uniform bound on its operator norm. In \cite{K14}, it was shown that under the uniform convexity assumption, the moment map indeed satisfies the uniform bound $||\Hess \varphi ||_{op} \leq \epsilon^{-1}$, and the conclusion follows. 
\end{proof}

It would also be possible to build a Stein kernel using the construction of \cite{Cha09, NPR10} and the optimal transport map sending the standard Gaussian measure onto $\mu$. Existence could be proved in the same setting, but there would be two main downsides: we do not have an analogue to Proposition \ref{est_lc_klartag} below for those maps, so we would only get a useful quantitative estimate in the uniformly convex setting, and due to the particular form of the construction of \cite{Cha09}, even in the latter setting the quantitative estimates would get worse. But we would still get existence of a Stein kernel that is bounded for uniformly log-concave measures. 

\section{Application to rates of convergence in the central limit theorem} \label{sect_clt}

We now show how the construction of Stein kernels discussed in the previous section leads to new estimates on the rate of convergence in the central limit theorem. The family of distances we shall consider to estimate the distance in the CLT are the Kantorovitch-Wasserstein distances from optimal transport theory, defined as
$$W_p(\mu, \nu) := \underset{\pi}{\inf} \hspace{1mm} \int{||x-y||_2^p\pi(dx, dy)}$$
where the infimum runs over all couplings of the measures $\mu$ and $\nu$. We refer to the textbook \cite{Vil03} for background on optimal transport. 

The following statement is a variant of a result of \cite{LNP15} on how $L^p$ bounds on a Stein kernel control Wasserstein distances to the standard Gaussian measure, which we shall prove in Section \ref{sect_bnd_wp}: 

\begin{proposition} \label{prop_bnd_wp}
Let $\tau$ be a Stein kernel for the probability measure $\mu$ on $\R^d$. Then for any $p \geq 2$ we have
$$W_p(\mu, \gamma) \leq C_p\left(\int{||\tau - \id||_{HS}^pd\mu}\right)^{1/p}$$
with $C_p = \left(\int{|x|^pd\gamma_1}\right)^{1/p}$ 
\end{proposition}

The original result of \cite{LNP15} bounds the Wasserstein distance of order $p$ by $\int{\underset{i,j}{\sum} \hspace{1mm} |\tau_{ij} - \delta_{ij}|^pd\mu}$, but with an extra prefactor that depends on the dimension. We shall use the above variant instead because it leads to a better dependence on the dimension for the quantitative CLT we shall later obtain. The prefactor $C_p$ behaves like $O(p)$.

These results mean that if we get estimates on $\Hess \varphi$, averaged out against $e^{-\varphi}$, we can deduce estimates on transport distances. It turns out that when $\mu$ is log-concave and compactly supported, such an estimate was already obtained by Klartag \cite{K14}: 

\begin{proposition} \label{est_lc_klartag} 
Let $\mu$ be a log-concave probability measure, supported on an open bounded convex set, and with a density bounded from above and below. Then the essentially-continuous convex function $\varphi$ for which $\mu$ is the moment measure is $C^2$ and satisfies for any $p \geq 1$ and any $\theta \in S^{d-1}$
$$\int{|\langle \Hess  \varphi(\nabla \varphi^*)\theta, \theta \rangle|^pd\mu} \leq 8^pp^{2p}\left(\int{(x \cdot \theta)^2d\mu}\right)^p.$$
\end{proposition}

We shall give a proof of Klartag's estimate in Section \ref{sect_proof_est_klartag}. It will be the same proof as in \cite{K14}, reformulated in a different language, which may be of interest to some readers. 

\begin{remark}
The results of \cite{K14} assume $C^{\infty}$ smoothness, relying on a result of \cite{BB13} to deduce $C^{\infty}$ smoothness of $\varphi$. Since we actually only need $C^2$ smoothness of $\varphi$, it turns out we only need continuity of the density. 
\end{remark}

Combining our construction of Stein kernels, Klartag's estimate and basic arguments from Stein's method, we get the following application to rates of convergence in the CLT: 

\begin{theorem} \label{clt_lp_rate_logconcave}
Let $\mu$ be an isotropic log-concave probability measure with strictly positive and continuous density on its support. Let $\mu_n$ be the law of $n^{-1/2}\underset{i=1}{\stackrel{n}{\sum}} \hspace{1mm} X_i$, where the $X_i$ are i.i.d. random variables with law $\mu$. Then for any $p \geq 2$ we have 
$$W_p(\mu_n, \gamma) \leq \tilde{C}(p)\frac{d}{n^{1/2}}$$
with $\tilde{C}(p)$ is a constant that depends only on $p$ (and which grows like $p^4$). In particular, this estimate does not depend on $\mu$. 
\end{theorem}

In the case $p =2$, the main result of \cite{CFP17} combined with the best currently-known estimate on the Poincar\'e constant of log-concave measures \cite{LV16} leads to a rate of convergence of the form $Cd^{3/4}n^{-1/2}$, which is better than the one we obtain here. The Kannan-Lovasz-Simonovits conjecture predicts that the Poincar\'e constant of isotropic log-concave measures is bounded by some universal constant, independently of the dimension, so we expect a rate of order $\sqrt{d/n}$. In a far more general setting, Bonis \cite{Bon16} proved the sharp rate of convergence for measures with moments of order $p+2$, and with a prefactor behaving like $d^{5/4}$ for general isotropic log-concave measures. In dimension one, Bobkov \cite{Bob17} also obtained the sharp rate for measures with finite moment. 

In the situation where $\mu$ is uniformly log-concave, the uniform estimate on the operator norm of our Stein kernel leads to an improved dependence on the dimension: 

\begin{theorem} \label{impro_clt_unif}
Let $\mu = e^{-V}$ be an isotropic probability measure with strictly positive and continuous density on its support, and assume that $\Hess V \geq \epsilon \operatorname{Id} $ for some $\epsilon > 0$. Then 
for any $p \geq 2$ we have 
$$W_p(\mu_n, \gamma) \leq C(p)\frac{\sqrt{d}}{\epsilon n^{1/2}}. $$
\end{theorem}

In this result, the dependence on the dimension is actually sharp, since it cannot be improved for product measures. It extends a result of \cite{CFP17} for $p=2$ to all $p \geq 2$. Once again, the constant $C(p)$ we obtain grows like $p^4$.

\begin{proof}[Proof of Theorem \ref{clt_lp_rate_logconcave}] 
We first work in the situation where $\mu$ has a compact support and a density bounded away from zero. Let $\tau = (\nabla^2 \varphi^*)^{-1}$, which we know is a Stein kernel for $\mu$. Then as is standard, $\tau_n(x) := \mathbb{E}\left[\frac{1}{n}\underset{k = 1}{\stackrel{n}{\sum}} \hspace{1mm} \tau(X_i) \left|\right. \frac{1}{\sqrt{n}}\underset{k = 1}{\stackrel{n}{\sum}} \hspace{1mm} X_i = x \right]$ is a Stein kernel for $\mu_n$. See for example the proof of Theorem 3.2 in \cite{CFP17} for a proof of this statement. Applying Jensen's inequality, we have
$$\int{|\tau_n - \id|_{HS}^p d\mu_n} \leq \int{\left|\frac{1}{n}\sum \tau(x_i) - \id\right|_{HS}^pd\mu^{\otimes n}(x_1,..,x_n)},$$
and Rosenthal's inequality for sums of independent random variables \cite{IS01} yields 
$$\int{|\tau_n - \id|_{HS}^p d\mu_n} \leq K_p^pn^{-p/2}\int{|\tau -\id|_{HS}^pd\mu}$$
with $K_p = O(p)$ the best constant in the Rosenthal inequality. This argument was already used in the discussion below Theorem 2.8 in \cite{LNP15}. We then have, given an orthonormal basis $(\theta_1,.., \theta_d)$ of $\R^d$,  
\begin{align*}
&\int{\left(\sum |\tau_{ij} - \delta_{ij}|^2\right)^{p/2}d\mu} \leq 2^p\left(d^{p/2} + \int{\left(\sum |\tau_{ij}|^2\right)^{p/2}d\mu}\right) \\
&\leq 2^p\left(d^{p/2} + \int{\left(\underset{i}{\sum} \hspace{1mm} \langle \tau \theta_i, \theta_i \rangle \right)^pd\mu}\right) \\
&\leq 2^p\left(d^{p/2} + d^{(p -1)}\sum \int{\langle \tau \theta_i, \theta_i\rangle^pd\mu}\right) \\
&\leq 2^pd^p(1 + 8^pp^{2p}). \\
\end{align*}
Hence
\begin{equation}
\int{|(\tau_n)_{ij} - \delta_{ij}|_{HS}^p d\mu_n} \leq K_pn^{-p/2}2^pd^{p}(1 + 8^pp^{2p}),
\end{equation} 
and therefore 
$$W_p(\mu_n, \gamma) \leq 2K_p(1 + 8^pp^{2p})^{1/p}C_p\frac{d}{n^{1/2}}.$$ 

For the general case, when the support of $\mu$ is not necessarily compact, we can take a sequence of compact sets $F_{\ell}$ that converge to the support of $\mu$, and apply our results to the restriction of $\mu$ to $F_{\ell}$ (renormalized to remain a centered, isotropic probability measure, so that $F_{\ell}$ has to be modified to take this into account, but this modification will remain convex and compact). The estimate on the Wasserstein distance does not depend on $F_{\ell}$, so that we can let $\ell$ go to infinity and the result remains valid. 
\end{proof}

The proof of Theorem \ref{impro_clt_unif} follows the exact same argument except that we use the improved bound of Corollary \ref{cor_est_stein_unif}, so we omit it. 

\subsection{Proof of Proposition \ref{est_lc_klartag}}

\label{sect_proof_est_klartag}

We shall now give a proof of Proposition \ref{est_lc_klartag}, omitting many computations taken from \cite{Kol14, KK17}. While it is not written in the same way as in \cite{K14}, it is the same proof, and we stress it is not due to us. We describe it in this form so that it is more easily readable for people with a knowledge of Bakry-Emery calculus

\begin{proof}[Proof of Proposition \ref{est_lc_klartag} ]
We introduce the Hessian metric on $\R^d$ given by the Riemannian metric tensor $g = \nabla^2 \varphi$. A result of Kolesnikov \cite{Kol14} asserts that when $\mu$ is log-concave and satisfies the regularity conditions of Theorem \ref{thm_moment_measures_smooth}, then the metric-measure space $M = (\R^d, g, e^{-\varphi})$ has Ricci curvature bounded from below by $1/2$. Moreover, if we consider the Laplacian on $M$, which is given by the formula 
$$L_{\varphi}f = \Tr(\nabla^2f (\nabla \varphi)^{-1}) + \nabla \log \rho (\nabla \varphi) \cdot \nabla f$$
then one can check that 
$$L_{\varphi}\partial_e \varphi = -\partial_e \varphi; \hspace{3mm} \Gamma(\partial_e \varphi) = \partial^2_{ee}\varphi; $$
where $\Gamma$ is the squared norm of the gradient with respect to the metric $g$. These computations can be found in \cite{Kol14, KK17}. We can then use tools from Bakry-\'Emery theory to obtain estimates on eigenfunctions of the Laplacian for spaces with positive curvature to deduce the desired bound. Indeed, if we denote by $P_t$ the semigroup acting on functions induced by $L_{\varphi}$, we have for any locally-lipschitz function $f$
$$\Gamma(P_t f) \leq \frac{1}{2(e^{t}-1)}P_t(f^2),$$
see Theorem 4.7.2 in \cite{BGL14}. Taking $f = \partial_e \varphi$, since it is an eigenfunction of $L_{\varphi}$ associated to the eigenvalue $1$, we have $P_t\partial_e \varphi = e^{-t}\partial_e\varphi$. Therefore $\Gamma (P_t f) = e^{-2t}\partial^2_{ee}\varphi$ and for any $t > 0$ and $p \geq 1$ we have
$$e^{-2pt}(\partial^2_{ee}\varphi)^p \leq \left(\frac{1}{2(e^{t}-1)}\right)^p(P_t((\partial_e \varphi)^2))^p \leq \left(\frac{1}{2(e^{t}-1)}\right)^pP_t((\partial_e \varphi)^{2p}).$$ Hence after integrating, for any $t > 0$ we have 
$$\int{(\partial^2_{ee}\varphi)^pe^{-\varphi}dx} \leq \left(\frac{e^{2t}}{2(e^{t}-1)}\right)^p\int{(\partial_e \varphi)^{2p}e^{-\varphi}dx}.$$
Taking $t = \ln 2$, the result then follows from the bound
$$||f||_{2p, e^{-\varphi}} \leq 2p||f||_{2, e^{-\varphi}}$$
for any eigenfunction of $L_{\varphi}$ associated with the eigenvalue $-1$, when a logarithmic Sobolev inequality with constant $1/2$ holds, see Section 5.3 of \cite{BGL14}. 
\end{proof}

\subsection{Proof of Proposition \ref{prop_bnd_wp}} \label{sect_bnd_wp}

Following the argument at the start of the proof of Proposition 3.4 in \cite{LNP15}, if we consider $X$ a random variable distributed according to $\mu$, $Z$ a standard Gaussian random variable independent of $X$, and $X_t := e^{-t}X + \sqrt{1 - e^{-2t}}Z$, we have when $p \geq 2$
\begin{equation*}
W_p(\mu, \gamma) \leq \int_0^{\infty}{\frac{e^{-2t}}{\sqrt{1 - e^{-2t}}}\mathbb{E}\left[\left(\mathbb{E}\left[\underset{i=1}{\stackrel{d}{\sum}}\left(\underset{j=1}{\stackrel{d}{\sum}}(\tau_{ij}(X) - \delta_{ij})Z_j\right)^2 \left|\right. X_t\right]\right)^{p/2}\right]^{1/p}dt}.
\end{equation*}
Applying Jensen's inequality, this is bounded by 
$$\int_0^{\infty}{\frac{e^{-2t}}{\sqrt{1 - e^{-2t}}}\mathbb{E}\left[\left(\underset{i=1}{\stackrel{d}{\sum}}\left(\underset{j=1}{\stackrel{d}{\sum}}(\tau_{ij}(X) - \delta_{ij})Z_j\right)^2\right)^{p/2}\right]^{1/p}dt}.$$
We can integrate in time, and we obtain
\begin{equation*}
W_p(\mu, \gamma) \leq \mathbb{E}[||(\tau - \id)Z||^p]^{1/p}
\end{equation*}
Since for a standard Gaussian random variable on $\R^d$ and a symmetric matrix $A$, we have $\mathbb{E}[||AZ||^p] \leq C_p^p||A||_{HS}^p$, this implies the desired bound. This moment bound can be proved by considering the case where $A$ is diagonal and applying the lower bound in the Marcinkiewicz-Zygmund inequality for sums of independent random variables with its sharp constant for symmetric random variables \cite{Fer14}. 

\section{A remark on the connection with weighted Poincar\'e inequalities}

The main result of \cite{CFP17} gives a construction of Stein kernels for measures satisfying a converse-weighted Poincar\'e, that is
$$\underset{c \in \R}{\inf} \hspace{1mm}\int{(f-c)^2\omega d\mu} \leq \int{|\nabla f|^2d\mu}$$
for some fixed weight function $w : \R^d \longrightarrow \R_+^*$, and all smooth scalar functions $f$. This notion generalizes the classical Poincar\'e inequality, which correpsonds to the case of constant weight function. 

In \cite{Sau18}, Saumard proved that in dimension one, the Stein kernel $\tau$ (if the density is sufficiently nice) gives rise to the following weighted Poincar\'e inequality
$$\Var_{\mu}(f) \leq \int{\tau (f')^2d\mu}.$$
Of course, this inequality exploits the fact that the (unique) Stein kernel in dimension one is nonnegative. Moreover Saumard showed that conversely the argument of \cite{CFP17} could be modified to prove that if such an inequality holds for some weight replacing $\tau$, then a Stein kernel exists (although the argument seems to guarantee the validity of the Stein identity for a possibly smaller class of test functions, depending on the behavior of the weight). The proof of the weighted Poincar\'e inequality makes use of the formula for Stein kernels in dimension one, and does not seem to readily extend to higher dimension. 

The Brascamp-Lieb inequality \cite{BL76} asserts that for any strictly log-concave probability measure $\nu = e^{-V}$, we have a weighted Poincar\'e-type inequality
$$\Var_{\nu}(f) \leq \int{\langle (\Hess V)^{-1}\nabla f, \nabla f \rangle d\nu}.$$
This inequality has many connections with geometric and functional inequalities (see for example \cite{BL00}), and has found applications in the study of long-time behavior of stochastic processes. 

It turns out that, combined with our construction of Stein kernels, this inequality immediately yields a multi-dimensional analogue of Saumard's result, at least as soon as the moment map is nice enough. Indeed, by taking $\varphi$ the moment map associated with a centered probability measure $\nu$, setting $\tau = \Hess \varphi(\nabla \varphi^*)$ as our Stein kernel, we have
\begin{align*}
\Var_{\nu}(f) &= \Var_{e^{-\varphi}}(f \circ \nabla \varphi) \\
&\leq \int{\langle (\nabla^2 \varphi)^{-1}\nabla^2\varphi \nabla f(\nabla \varphi), \nabla^2\varphi \nabla f(\nabla \varphi)\rangle e^{-\varphi}dx} \\
&= \int{\langle \nabla f(\nabla \varphi), \nabla^2\varphi \nabla f(\nabla \varphi)\rangle e^{-\varphi}dx} \\
&= \int{\langle \tau \nabla f, \nabla f\rangle d\nu}. 
\end{align*}

Hence we obtain a multi-dimensional generalization of Saumard's result, albeit for a specific choice of Stein kernel. Note that it is not clear at all that we should expect general Stein kernels to take values that are positive matrices, so a similar weighted Poincar\'e inequality for any kernel might not be true. In particular, it is not clear at all that the Stein kernels constructed in \cite{CFP17} have such a property. But in dimension one, since all possible Stein kernels match when the measure has a density, this in particular gives an alternative proof of Saumard's result.

\section{Transporting Stein kernels for other reference measures}

The abstract setup of Stein's method can be generalized to cover non-gaussian reference measures \cite{Bar90}. If we wish to compare some measure $\mu$ to a reference measure $\mu_0 = e^{-V}dx$, say for a smooth function $V$ that is finite everywhere, then $\mu_0$ is characterized by the integration by parts formula
$$\int{\nabla V \cdot f d\mu_0} = \int{\Tr(\nabla f)d\mu_0}$$
which leads to a definition of Stein kernel as a matrix-valued function such that
\begin{equation}
\int{\nabla V \cdot f d\mu} = \int{\langle \tau, \nabla f \rangle d\mu}. 
\end{equation}

Assume that $V$ is convex, $C^2$ and that $\Hess V > 0$, and let $\mu_V$ be the pushforward of $\mu$ by $\nabla V$. Then for any vector-valued smooth function $f$, defining $g(x) = f(\nabla V^*(x))$, we have
\begin{align*}
\int{\nabla V(x) \cdot f(x) d\mu} &= \int{\nabla V(x) \cdot g(\nabla V(x)) d\mu} \\
&= \int{x \cdot g(x) d\mu_V}, 
\end{align*}
so if we take $\tilde{\tau}_{V, \gamma}$ to be a Stein kernel for $\mu_V$ \emph{with respect to the gaussian measure} (assuming for now it exists), we get
\begin{align*}
\int{\nabla V(x) \cdot f(x) d\mu} &= \int{\langle \tilde{\tau}_{V, \gamma}, \nabla g \rangle d\mu_V} \\
&= \int{\langle \tilde{\tau}_{V, \gamma}(\nabla V(x)), \nabla g(\nabla V(x)) \rangle d\mu} \\
&= \int{\langle \tilde{\tau}_{V, \gamma}(\nabla V(x)), (\Hess V(x))^{-1}\nabla f(x) \rangle d\mu} \\
&= \int{\langle \tilde{\tau}_{V, \gamma}(\nabla V(x))(\Hess V(x))^{-1}, \nabla f(x) \rangle d\mu}
\end{align*}
and therefore $\tilde{\tau}_{V, \gamma}(\nabla V(x))(\Hess V(x))^{-1}$ is a Stein kernel for $\mu$ relative to $\mu_0$. To be valid, in addition to the regularity and convexity assumptions on $V$, this arguments requires that $\tilde{\tau}_{V, \gamma}$ exists. It is okay if it is only defined in the sense of distributions (since $\nabla V$ is smooth and bijective, composing a distribution with it is possible). 

\vspace{3mm}

\underline{\textbf{ Acknowledgments }} : This work was partly made while the author was in residence at the Institut Henri Poincar\'e in July 2017, and at the Mathematical Sciences Research Institute in Berkeley for part of the fall 2017 semester. The author benefited from support from the France Berkeley Fund, the ANR Project ANR-17-CE40-0030 - EFI, and ANR-11-LABX-0040-CIMI within the program ANR-11-IDEX-0002-02. I thank Dimitri Shlyakhtenko for the discussion that initiated this work, and Thomas Bonis for pointing out a way to improve the dependence on the dimension of the results. I also thank Dario Cordero-Erausquin, Thomas Courtade, Bo'az Klartag, Michel Ledoux, Emanuel Milman and the anonymous referee for their advice and comments, and Adrien Saumard for sharing a preliminary version of \cite{Sau18}.

\end{document}